\documentclass{gtart}

\def\ifplaintex{\expandafter\ifx\csname documentclass\endcsname\relax}

\expandafter\ifx\csname epsfbox\endcsname\relax\input epsf\fi

\ifplaintex 
\else
\fi

\def\gt{{\mathsurround=0pt\it $\cal G\mskip-2mu$eometry \&\ 
$\cal T\!\!$opology}}        

\def\gtp{{\mathsurround=0pt\it $\cal G\mskip-2mu$eometry \&\ 
$\cal T\!\!$opology $\cal P\!$ublications}}  

\def\lognumber#1{\def\thelognumber{#1}}
\def\volumenumber#1{\def\thevolumenumber{#1}}
\def\papernumber#1{\def\thepapernumber{#1}}
\def\volumeyear#1{\def\thevolumeyear{#1}}

\def\pagenumbers#1#2{\def\startpage{#1}\def\finishpage{#2}}
\def\published#1{\def\publishdate{#1}}
\def\proposed#1{\def\theproposer{#1}}
\def\seconded#1{\def\theseconders{#1}}
\def\received#1{\def\receiveddate{#1}}
\def\revised#1{\def\reviseddate{#1}}
\def\accepted#1{\def\accepteddate{#1}}
\def\asciititle#1{\def\theasciititle{#1}}

\def\coverauthors#1{\def\thecoverauthors{#1}}
\def\asciiauthors#1{\def\theasciiauthors{#1}}
\def\asciiaddress#1{\def\theasciiaddress{#1}}

\long\def\asciiabstract#1{\long\def\theasciiabstract{#1}}
\def\asciikeywords#1{\def\theasciikeywords{#1}}

\let\\\par\let\thevolumenumber\relax\let\thepapernumber\relax
\let\thevolumeyear\relax\let\thesamplenumber\relax\let\startpage\relax
\let\finishpage\relax\let\publishdate\relax\let\receiveddate\relax
\let\reviseddate\relax\let\accepteddate\relax\let\theasciititle\relax
\let\theasciiauthors\relax\let\theasciiaddress\relax
\let\theasciiabstract\relax\let\theasciikeywords\relax
\let\theasciiemail\relax\let\theshortauthors\relax\let\theshorttitle\relax
\let\thecoverauthors\relax

\long\def\maketitlep{   

\count0=\startpage

\gt\hfill      
\hbox to 77pt{\vbox to 0pt{\vglue -15pt\epsfbox{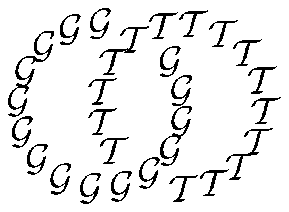}\vss}\hss}
\break
{\small\ifx\thesamplenumber\relax 
Volume \else Sample
\fi\thevolumenumber\ (\thevolumeyear)
\startpage--\finishpage\nl
Published: \publishdate}
\vglue 0.5truein plus 0.4fil minus 0.1truein

{\parskip=0pt\leftskip 0pt plus 1fil\def\\{\par\smallskip}{\ifplaintex\large
\else\Large\fi\bf\thetitle}\par\medskip}   

\vglue 0pt plus 0.1fil 

{\parskip=0pt\leftskip 0pt plus 1fil\def\\{\par}{\sc\theauthors}
\par\medskip}

\vglue 0pt plus 0.1fil 

{\small\parskip=0pt\let\newline\\
{\leftskip 0pt plus 1fil\def\\{\par}{\sl\theaddress}\par}
\expandafter\ifx\theemail\relax    
\relax\else\vglue 5pt plus 0.02fil minus 2pt\def\\{\stdspace{\rm 
and}\stdspace} 
\cl{Email:\stdspace\tt\theemail}\fi
\ifx\theurl\relax                  
\relax\else\vglue 5pt plus 0.02fil minus 2pt\def\\{\stdspace{\rm 
and}\stdspace}
\cl{URL:\stdspace\tt\theurl}\fi\par}

\vglue 7pt plus 0.3fil minus 3pt

{\bf Abstract}
\vglue 5pt plus 0.1fil minus 2pt

\theabstract

\vglue 7pt plus 0.3fil minus 3pt

{\bf AMS Classification numbers}\quad Primary:\quad \theprimaryclass

Secondary:\quad \thesecondaryclass

\vglue 5pt plus 0.3fil minus 2pt

{\bf Keywords:}\quad \thekeywords

\vglue 10pt plus 0.5fil minus 5pt

{\small  Proposed: \theproposer\hfill Received: \receiveddate\nl
Seconded: \theseconders\hfill 
\ifx\reviseddate\relax                         
Accepted: \accepteddate                        
\else
Revised: \reviseddate                          
\fi}
\eject
}       

\let\maketitlepage\maketitlep
\let\maketitle\maketitlepage

\font\phead=cmsl9 scaled 950
\font=cmsl9 scaled 1050
\font\pnum=cmbx10 scaled 913
\font\lnum=cmbx10 
\font\pfoot=cmsl9 scaled 950
\font=cmsl9 scaled 1050
\ifplaintex
\headline{\vbox to 0pt{\vskip -4.5mm\line{\small\phead\ifnum
\count0=\startpage ISSN 1364-0380 (on line)
1465-3060 (printed) \hfill {\pnum\folio}\else\ifodd\count0\def\\{ }%
\ifx\theshorttitle\relax\thetitle\else\theshorttitle\fi\hfill{\pnum\folio}
\else\def\\{ and }{\pnum\folio}\hfill\ifx\theshortauthors\relax\theauthors
\else\theshortauthors\fi\fi\fi}\vss}}
\footline{\vbox to 0pt{\vglue 0mm\line{\small\pfoot\ifnum\count0=\startpage
\copyright\ \gtp\hfill\else
\gt, Volume \thevolumenumber\ (\thevolumeyear)\hfill\fi}\vss
}}
\else
\makeatletter
\def\@oddhead{{\small\lhead\ifnum\count0=\startpage ISSN 1364-0380 (on line)
1465-3060 (printed) \hfill {\lnum\number\count0}\else\ifodd\count0
\def\\{ }\ifx\theshorttitle\relax \thetitle \else\theshorttitle\fi\hfill
{\lnum\number\count0}\else\def\\{ and }{\lnum\number\count0}
\hfill\ifx\theshortauthors\relax 
\theauthors\else\theshortauthors\fi\fi\fi}}\def\@evenhead{@oddhead}
\def\@oddfoot{\small\lfoot\ifnum\count0=\startpage\copyright\ \gtp\hfill\else
\gt, Volume \thevolumenumber\ (\thevolumeyear)\hfill\fi}
\def\@evenfoot{@oddfoot}
\makeatother
\fi

\newwrite\gtoutfile
\long\gdef\makeheadfile{  
{\def\\{, }\def\s{ }
\immediate\openout\gtoutfile head.xxx
\immediate\write\gtoutfile{To: math@arxiv.org}
\immediate\write\gtoutfile{Subject: put OR rep NNNNN:pppp}
\immediate\write\gtoutfile{--text follows this line--}
\immediate\write\gtoutfile{Proxy-for: \ifx\theasciiauthors\relax
\theauthors\else\theasciiauthors\fi\s<\ifx\theasciiemail\relax\theemail\else\theasciiemail\fi>}
\immediate\write\gtoutfile{\noexpand\\}
\immediate\write\gtoutfile{Authors: \ifx\theasciiauthors\relax
\theauthors\else\theasciiauthors\fi}
{\def\\{ }\immediate\write\gtoutfile{Title: \ifx\theasciititle\relax
\thetitle\else\theasciititle\fi}}
\immediate\write\gtoutfile{Subj-class: GT or GR or SG or ...}
\immediate\write\gtoutfile{MSC-class: \theprimaryclass\ifx\thesecondaryclass\relax\else, \thesecondaryclass\fi}
\immediate\write\gtoutfile{Journal-ref: Geom. Topol. \thevolumenumber\s
(\thevolumeyear) \startpage-\finishpage}
\immediate\write\gtoutfile{Comments: Published in Geometry and Topology at}
\immediate\write\gtoutfile{    http://www.maths.warwick.ac.uk/gt/GTVol\thevolumenumber/paper\thepapernumber.abs.html}
\immediate\write\gtoutfile{\noexpand\\}
\immediate\write\gtoutfile{}
\ifx\theasciiabstract\relax
\immediate\write\gtoutfile{\theabstract}\else
\immediate\write\gtoutfile{\theasciiabstract}\fi
\immediate\write\gtoutfile{}
\immediate\write\gtoutfile{\noexpand\\}
\immediate\write\gtoutfile{}
\immediate\closeout\gtoutfile}}  

\def\maketitlepage{\maketitlep\makeheadfile}
\let\maketitle\maketitlepage

\lognumber{208}
\volumenumber{6}\papernumber{7}\volumeyear{2002}
\pagenumbers{195}{218}
\received{1 October 2001}
\revised{11 March 2002}
\published{27 April 2002}
\accepted{26 April 2002}

\proposed{Yasha Eliashberg}
\seconded{Tomasz Mrowka, John Morgan}

\usepackage{amsmath, amssymb}
\usepackage[all]{xy}

\newcommand{\R}{{\mathbb R}}

\newcommand{\p}{{\partial}}

\newcommand{\Jj}{{\mathcal J}}

\newcommand{\Nn}{{\mathcal N}}

\newcommand{\Q}{{\mathbb Q}}
\newcommand{\N}{{\mathbb N}}
\newcommand{\Pp}{{\bf P}}

\newcommand{\Z}{{\mathbb{Z}}}
\newcommand{\C}{{\mathbb C}}

\newcommand{\al}{{\alpha}}

\newcommand{\be}{{\beta}}

\newcommand{\Om}{{\Omega}}
\newcommand{\om}{{\omega}}

\newcommand{\ga}{{\gamma}}

\newcommand{\la}{{\lambda}}
\newcommand{\si}{{\sigma}}
\newcommand{\La}{{\Lambda}}

\newcommand{\ph}{{\varphi}}

\DeclareMathOperator{\Diff}{Diff}
\DeclareMathOperator{\Map}{Map}

\newtheorem{theorem}{Theorem}[section]

\newtheorem{corollary}[theorem]{Corollary}
\newtheorem{lemma}[theorem]{Lemma}

\newtheorem{proposition}[theorem]{Proposition}
\theoremstyle{remark}
\newtheorem{remark}[theorem]{Remark}
\newtheorem{definition}[theorem]{Definition}

\begin{document}

\title{Homotopy type of symplectomorphism groups\\of $S^2 \times
S^2$}
\asciititle{Homotopy type of symplectomorphism groups of S^2 X S^2}

\author{S\'{\i}lvia  Anjos}
\coverauthors{S\noexpand\'{\noexpand\i}lvia  Anjos}
\address{ Departamento de Matem\'atica\\Instituto Superior T\'{e}cnico, Lisbon, Portugal}
\asciiauthors{Silvia  Anjos}
\asciiaddress{Departamento de Matematica\\Instituto 
Superior Tecnico, Lisbon, Portugal}
\email{sanjos@math.ist.utl.pt}

\begin{abstract}

In this paper we discuss the topology of the symplectomorphism group
of a product of two 2--dimensional spheres when the ratio of their
areas lies in the interval (1,2]. More precisely we compute the
homotopy type of this symplectomorphism group and we also show that
the group contains two finite dimensional Lie groups generating the
homotopy. A key step in this work is to calculate the mod 2 homology
of the group of symplectomorphisms. Although this homology has a
finite number of generators with respect to the Pontryagin product, it
is unexpected large containing in particular a free noncommutative
ring with 3 generators.
\end{abstract} 

\asciiabstract{In this paper we discuss the topology of the
symplectomorphism group of a product of two 2-dimensional spheres
when the ratio of their areas lies in the interval (1,2]. More
precisely we compute the homotopy type of this symplectomorphism group
and we also show that the group contains two finite dimensional Lie
groups generating the homotopy. A key step in this work is to
calculate the mod 2 homology of the group of
symplectomorphisms. Although this homology has a finite number
of generators with respect to the Pontryagin product, it is unexpected
large containing in particular a free noncommutative ring with 3
generators.}

\asciikeywords{Symplectomorphism group, Pontryagin ring, homotopy 
equivalence}

\maketitlepage

 \section{Introduction} In general symplectomorphism groups
are thought to be intermediate objects between Lie groups and full
groups of diffeomorphisms. Although  very little is known about
the topology of groups of diffeomorphisms, there are some cases when the corresponding symplectomorphism groups are more
understandable. For example, nothing is known about the group of
compactly supported diffeomorphisms of $\R^4$, but in 1985, Gromov
showed in  \cite{G} that the group of compactly supported
symplectomorphisms of $\R^4$ with its standard symplectic
structure is contractible. He also showed that the
symplectomorphism group of a product of two 2--dimensional spheres
that have the same area has the homotopy type of a Lie group.

 More
precisely, let $M_{\la}$ be the symplectic manifold $(S^2 \times
S^2, \om_{\la}=(1+\la) \si_0 \oplus \si_0)$ where $0 \leq \la \in
\R$ and $\si_0$ is the standard area form on $S^2$ with total area
equal to 1. Denote by $G_{\la}$ the group of symplectomorphisms of
$M_{\la}$ that act as the identity on $H_2(S^2 \times S^2; \Z)$.
Gromov proved that $G_0$ is connected and it is homotopy
equivalent to its subgroup of standard isometries $SO(3) \times
SO(3)$. He also showed that this would no longer hold when one
sphere is larger than the other, and in \cite{M1} McDuff constructed
explicitly an element of infinite order in $H_1(G_{\la})$, $\la
>0$. The main tool in their proofs is to look at the action of
$G_{\la}$ on the contractible space $\Jj_{\la}$ of
$\omega_{\la}$--compatible almost complex structures.

 Abreu and
McDuff in \cite{AM} calculated the rational cohomology of these
symplectomorphism groups and confirmed that these groups could not
be homotopic to Lie groups. In particular they computed the
cohomology algebra $H^{\ast}(G_{\la},\Q)$ for every $\la$. For
each integer $\ell \geq 1$ we have $$ H^{\ast}(G_{\la};\Q) = \La
(t,x,y) \otimes \Q[w_\ell], \mbox{ when $\ell-1 < \la \leq \ell$
}$$ where $\La (t,x,y)$ is an exterior algebra over $\Q$ with
generators $t$ of degree 1, and $x,y$ of degree 3 and $\Q[w_\ell]$
is the polynomial algebra on a generator $w_\ell$ of degree
$4\ell$ that is made from $x,y$, and $t$ via higher Whitehead
products. The generator $w_\ell$ is fragile, in the sense that it
disappears (ie, becomes null cohomologous) when $\la$ increases.
Moreover they showed that the rational homotopy type of $G_{\la}$
changes precisely when the ratio of the size of the larger to the
smaller sphere passes an integer.

In this paper we show that when $0 < \la \leq 1$  the whole
homotopy type  of $G_{\la}$ (rather than just its rational part)
is generated  by
 its subgroup of isometries $SO(3) \times SO(3)$ and by this new element of
infinite
 order constructed by McDuff.  More precisely we will calculate the homotopy type of
  $G_{\la}$:
\begin{theorem}
\label{tm:homtype}
If $0 < \la \leq 1$, $G_{\la}$ is homotopy equivalent to the product
$X= L \times S^1 \times  SO(3) \times SO(3)$ where $L$ is the loop space of
the
suspension of the smash  product $S^1 \wedge SO(3)$.
\end{theorem}
In this product of $H$--spaces\footnote{$X$ is an $H$--space if there is a map $\mu \co X \times X \rightarrow X$ such that  $\mu \circ i_1 \simeq 1$ and  $\mu \circ i_2 \simeq 1$ where $i_1$ and $i_2$ are the inclusions $i_1(x)=(x, \ast)$ and $i_2(x)=(\ast, x)$,  $\cong$ means homotopy equivalent and $* \in X$ is a base point.} one of the $SO(3)$ factors
corresponds to rotation in one of the spheres, the other
represents the diagonal in $SO(3) \times SO(3)$, and the $S^1$
factor corresponds to the generator in
 $H_1(G_{\la})$ described by Gromov and McDuff.  This new
element of infinite order represents
  a $S^1$--action that commutes with the diagonal action of $SO(3)$, but not
with
  rotations in each one of the spheres. The loop space
  $L=\Omega\Sigma (S^1 \wedge SO(3))$ appears as the result of that
non-commutativity.

Although this space $X$ is an $H$--space, its multiplication is not
the same as on $G_{\la}$. This can be seen by comparing the
Pontryagin products on integral homology.

The main steps in the proof of this theorem determine the
organization of the paper. Therefore in Section 2 we have the first
main result which is the calculation of the  mod 2 homology ring
$H_{\ast} (G_{\la} ; \Z_{2} )$. Recall that the product structure
in $H_{\ast} (G_{\la} ; \Z_{2} )$, called Pontryagin product, is
induced by the product in $G_{\la}$. Denote by $\La (y_1,...,y_n)$
the exterior algebra over $\Z_2$ with generators $y_i$ where this
means that $y_i^2=0$ and $y_iy_j=y_jy_i$ for all $i,j$, and by
$\Z_2 \langle x_1,...,x_n \rangle$ the free noncommutative algebra
over $\Z_2$ with generators $x_j$. Recall that
$H_{\ast}(SO(3);\Z_2) \cong \La (y_1,y_2)$.

\begin{theorem}\label{tm:algiso}
If $0 <\la \leq 1$ then there is an algebra isomorphism $$H_{\ast}
(G_{\la} ; \Z_{2} )= \La(y_1,y_2) \otimes \Z_2\langle
t,x_1,x_2\rangle /R$$ where $\deg y_i= \deg x_i =i$, $\deg t=1$ and $R$ is the
set of relations $ \{ t^2=x_i^2=0, x_1x_2=x_2x_1 \}$.
\end{theorem}

\noindent The notation implies that $y_i$ commutes with $t$ and
$x_i$. We see that  $H_{\ast} (G_{\la} ; \Z_{2} )$ contains $\La
(x_1,x_2)$ which appears from rotation in the first sphere, $\La
(y_1,y_2)$ which represents  the diagonal in $SO(3) \times SO(3)$
plus the new generator  in $H_1(G_{\la})$, $\la>0$, that we denote by $t$.  From the inclusion of
the subgroup of isometries $SO(3) \times SO(3)$ in $G_{\la}$ we
have classes $x_1,x_2,x_3=x_1x_2 \in H_{\ast} (G_{\la} ; \Z_{2} )$ in
dimensions 1,2 and 3 respectively, representing the rotation in
the first
 factor.  The new generator $t$ in $H_1(G_{\la})$
does not
 commute with $x_i$,
therefore  we have a nonzero class defined as the commutator and
represented by $x_it+tx_i$ for $i=1,2,3$. It is easy to understand
what these  classes are in homotopy.  For example, $x_1$ is a
spherical class, so it represents an element in $\pi_1 (G_{\la})$
and $x_1t+tx_1$ corresponds to  the Samelson product $[t,x_1]\in
\pi_2(G_{\la}) $. This is given by the map $$S^2=S^1 \times S^1 /
S^1 \vee S^1 \rightarrow G_{\la} $$ induced by the commutator
$$S^1 \times S^1 \rightarrow G_{\la} \co  (s,u) \mapsto
t(s)x_1(u)t(s)^{-1}x_1(u)^{-1}.$$
 Although the mod 2 homology has a
finite number of generators with respect to the Pontryagin product
we will see it is very large containing in particular
a free noncommutative ring on 3 generators, namely the commutators $x_it+tx_i$, $i=1,2,3.$

The proof of the theorem generalizes Abreu's work and is based on
the fact, proved by Abreu in \cite{A}, that the space $\Jj_{\la}$, of almost complex structures
on $S^2 \times S^2$ compatible with $\om_{\la}$, is a stratified
space with two strata $U_0$ and $U_1$, where $U_0$ is the open
subset of $\Jj_{\la}$ consisting of all $J \in \Jj_{\la}$  for
which the homology classes  $E=[S^2 \times \{ pt \}]$ and $F=[\{
pt \} \times S^2]$ are both represented by $J$--holomorphic spheres
and  its complement $U_1$ is a submanifold of codimension 2. More
precisely, $U_1$ consists of all $J \in \Jj_{\la}$ for which the
homology class of the antidiagonal  $E-F$ is represented by a
$J$--holomorphic sphere.

In Section 3 we start by giving some considerations about torsion in
$H_{\ast}(G_{\la};\Z)$. In particular we establish that
$H_{\ast}(G_{\la};\Z)$ has only 2--torsion. Then we define a map $f$
between $G_{\la}$ and the product $X=L \times S¹ \times SO(3)
\times SO(3)$ and prove it is in fact an homotopy equivalence. This is
obtained from three topological facts: (i) it is enough to find a
$\Z$--homology isomorphism from another $H$--space; (ii) a
$\Z$--homology isomorphism is implied by an isomorphism with all field
coefficients; (iii) homology is computed via Leray--Hirsch for two
fibrations of $G_{\la}$ over $U_0$ and $U_1$ and a model space is
built using universal properties for maps from loop spaces to
topological monoids.

\rk{Acknowledgements} The results of this paper are part of my
doctoral thesis, written at State University of New York at Stony
Brook under the supervision of D McDuff. I would like to thank her for
all her guidance, advice and support.  I also want to thank M Abreu
for useful corrections on a version of this paper.

The author acknowledges support from FCT (Funda\c{c}\~{a}o para a
Ci\^encia e Tecnologia) -- Programa Praxis, Funda\c{c}\~{a}o
Luso-Americana para o Desenvolvimento and Funda\c{c}\~{a}o Calouste
Gulbenkian.

 \section{The Pontryagin ring $H_{\ast}(G_{\la};\Z_2) $ }
 \label{sec:ring} 

Recall that for any group $G$ the product $\phi \co G \times G
\rightarrow G$ induces a product in homology $$ H_{\ast}(G;\Z_2 )
\otimes H_{\ast}(G;\Z_2 ) \xrightarrow{ \times} H_{\ast}(G \times
G;\Z_2 ) \xrightarrow{\phi_{\ast}} H_{\ast}(G;\Z_2 ) $$ called the
Pontryagin product, that we will denote by ``.''. Every time it is
clear from the context we will suppress this for simplicity of
notation. In this section we will compute the ring structure on
$H_{\ast}(G_{\la};\Z_2) $ induced by this product. Unless noted
otherwise we assume $\Z_2$ coefficients throughout.

\subsubsection{Geometric description}

As we mentioned in the introduction, Abreu proved in \cite{A} that if
$ 0 < \la \leq 1 $ the space of almost complex structures compatible
with $\om_{\la}$, $\Jj_{\la}$, is a stratified space with two strata
$U_0$ and $U_1$, where $U_0$ is open and dense and $U_1$ has
codimension 2. $U_i$ is the set consisting of all $J \in \Jj_{\la}$
for which the class $E-iF$ is represented by a $J$-holomorphic sphere,
where $E$ denotes the homology class of $S^2 \times \{pt\}$ and $F$
denotes the fiber class $\{pt\} \times S^2$. The group of
symplectomorphisms $G_{\la}$ acts on $\Jj_{\la}$ by
conjugation. Moreover the group $G_{\la}$ has finite dimensional
subgroups $K_i$, with $i=0,1$, acting on $M_{\la}$, where $K_0=SO(3)
\times SO(3)$ corresponds to the standard K\"{a}hler action of $SO(3)
\times SO(3)$ on $S^2 \times S^2$ with complex structure the standard
split structure $J_0= j_0 \oplus j_0$ and $K_1=SO(3) \times S^1$ is a
K\"{a}hler action for a complex structure $J_1 \in U_1$ with the
property that the unique $J_1$--holomorphic representative $C_2$ for
the class $E-F$ is fixed by the $S^1$ part of the action (see
below). The $SO(3)$ part of this action is the same as the diagonal
$SO(3)$ action on $S^2 \times S^2$.

The next step is to identify each stratum $U_i$ of $\Jj_{\la}$ with the quotient of $G_{\la}$ by the isometry group $K_i$. The result was proved by Abreu in \cite{A} and is the following:

\begin{proposition}   The stratum $U_i \in \Jj_{\la}$  is weakly homotopy equivalent to the quotient $G_{\la}/K_i$, $i=0$ or 1.\end{proposition}

Now we can give a brief geometric description of the $S^1$ part of the action in $K_1$ corresponding to the element of infinite order in $\pi_1(G_{\la})$ constructed by McDuff in \cite{M1}. The complex structure $J_1$ is tamed by $\om_{\la}$ and the complex manifold $ (S^2
\times S^2, J_1) $ is biholomorphic to the projectivization $\Pp ( \mathcal{O} (2) \oplus \C )$ over $S^2$.
Here $\mathcal{O} (2) $ is a complex line bundle over $S^2$ with
first Chern class 2. This bundle has two natural sections, $\Pp (
\{ 0\}  \oplus \C )$ and $\Pp ( \mathcal{O} (2) \oplus \{ 0\} )$, which represent the classes    $E+F$ (the diagonal in $S^2 \times S^2$) and $E-F$ (the antidiagonal in $S^2 \times S^2$). The element of infinite order in
$\pi_1(G_{\la})$ acts on this fibration  by rotation on the fibers
and leaving fixed the sections corresponding to the classes of the
diagonal and antidiagonal. We see that this element is in the
stabilizer of $J_1$ in $G_{\la}$, because this rotation is a
complex operation. Moreover for each $J \in U_0$ in a
neighborhood of $U_1$ the action of $t \in \pi_1(G_{\la})$ on $J$
gives a loop around $U_1$ which represents the link of $U_1$ in
$U_0$.

\subsection{Relation between $H_{\ast}(G_{\la})$ and
$H_{\ast}(U_i)$: additive version}

The fact that  $U_1$
is a codimension 2 submanifold of $\Jj_{\la}$ implies  that
$U_0=\Jj_{\la} - U_1$ is connected. This means that $G_{\la}$ is connected, which in turn implies that $U_1$ is also connected. Hence \begin{displaymath}
H_{0}(U_0;{\Z}_{2}) \cong {\Z}_2 \cong H_{0}(U_1;\Z_{2}).
\end{displaymath} {\rm Just as M.Abreu showed in \cite{A} we
still have for $p\geq 1$, \begin{equation} H_{p}(U_0;{\Z}_{2})
\cong H_{p-1}(U_1;{\Z}_{2}). \label{eq:MV} \end{equation} This
already implies that $H_{1}(U_0;\Z_{2}) \cong \Z_{2}$.} \noindent Now  consider the following principal fibrations \begin{equation}
\label{eq:diag}
 \xymatrix{ {K_0}
\ar[r]^{i_0} & {G_{\la}} \ar[d]^{p_0} \\
 & U_0 }   \hspace{1in}
\xymatrix{ {K_1} \ar[r]^{i_1} & {G_{\la}} \ar[d]^{p_1} \\
 & U_1 }
\end{equation}
 where $K_i$ is the identity component of the stabilizer of $J_i$ in $G_{\la}$. As we stated before $K_0$ is the subgroup
$SO(3) \times SO(3)$ and $K_1$ is isomorphic to $S^1 \times
SO(3)$.

The following proposition was proved by Abreu for rational
coefficients but we need it for $\Z_2$ coefficients.
\begin{proposition} \label{prop:map} Let ${\Diff}_0 (S^2 \times
S^2)$ denote the group of diffeomorphisms of $S^2 \times S^2$ that
act as the identity on $H_2 (S^2 \times S^2, \Z )$. The inclusion
$$ i\co K_0=SO(3) \times SO(3) \longrightarrow \Diff_0 (S^2 \times
S^2) $$ is injective in homology.
\end{proposition} \begin{proof} As in \cite{A} we define a map $$
F\co {\Diff}_0 (S^2 \times S^2) \longrightarrow \Map_1(S^2) \times
\Map_1 (S^2) $$  where $\Map_1(S^2)$  is the space of all
orientation preserving self-homotopy equivalences of $S^2$. Given
$\varphi \in \Diff_0 (S^2 \times S^2)$ we define a self map of
$S^2$, denoted by $ \tilde{\varphi}_1$, via the composite
$$\tilde{\varphi}_1 \co S^2 \stackrel{i_1}{\rightarrow} S^2 \times
S^2 \stackrel{\varphi}{ \rightarrow} S^2 \times S^2
\stackrel{\pi_1}{\rightarrow} S^2,$$ where $i_1$, respectively
$\pi_1$, denote inclusion into, respectively projection onto, the
first $S^2$ factor of $S^2 \times S^2$. Because $\varphi$ acts as
the identity on $H_2 (S^2 \times S^2, \Z )$, $ \tilde{\varphi}_1$
is an orientation preserving self homotopy equivalence of $S^2$,
ie, $ \tilde{\varphi}_1 \in \Map_1 (S^2) $. Defining $
\tilde{\varphi}_2$ in an analogous way using the second $S^2$ factor
of $S^2 \times S^2$, we have thus constructed the desired map
given by $$ \varphi \mapsto \tilde{\varphi}_1 \times
\tilde{\varphi}_2.$$ It is clear from the construction that $F$
restricted to $SO(3) \times SO(3)$ is just the inclusion $$ SO(3)
\times SO(3) \longrightarrow  \Map_1(S^2) \times \Map_1 (S^2) $$
Now we use the following theorem (see \cite{H}).
\begin{theorem} \label{thm:heq}The space of orientation preserving
self-homotopy equivalences on the 2--sphere has  the homotopy type
of $SO(3) \times \Om$, where $\Om= \Tilde{\Om}^2_0 (S^2)$ is the
universal covering space for the component in the double loop
space on $S^2$ containing the constant based map. \end{theorem}

This proves that $SO(3)$ is not homotopy equivalent to
$\Map_1(S^2)$ but we  have, using the K\"{u}nneth formula with
field coefficients, $$ H_{\ast}(SO(3) \times \Om) \cong
H_{\ast}(SO(3)) \otimes H_{\ast}(\Om)\cong H_{\ast}(\Map_1(S^2))
$$ thus the map $$ i_{\ast}\co H_{\ast}(SO(3)) \longrightarrow
H_{\ast}(\Map_1(S^2))$$ induced by injection is injective  for any
field coefficients. \end{proof}

It is proved by D McDuff in \cite{M1} that the generator of the $\Z$
factor in
 $\pi_1(G_{\la})$ lies in $\pi_1(K_1)$. This means that the generator of the
  $S^1$--action in  $\pi_1(K_1)$ maps to a generator of infinite order in
  $\pi_1(G_{\la})$. Thus the map $${i_1}_{\ast}\co H_{\ast}(K_1 )
\longrightarrow H_{\ast}(G_{\la})  $$ induced by inclusion is
injective. Since we are working over a field, the cohomology is
the dual of homology, thus from the above and Proposition \ref{prop:map}
the maps $$i^{\ast}_0 \co H^{\ast}(G_{\la}) \longrightarrow H^{\ast}(K_0)$$ 
and $$i^{\ast}_1 \co H^{\ast}(G_{\la})
\longrightarrow H^{\ast}(K_1)$$ induced by inclusions $i_0$
and $i_1$ are surjective.

From the Leray--Hirsch Theorem it follows that the spectral sequences of the
fibrations
collapse at the $E_2$--term, and we have the following vector space
isomorphisms \begin{equation} \label{eq:cxa} H^{\ast}(G_{\la} )
\cong H^{\ast}(U_0)  \otimes  H^{\ast}(K_0) \end{equation}
\begin{equation} \label{eq:ca} H^{\ast}(G_{\la}) \cong
H^{\ast}(U_1 ) \otimes H^{\ast}(K_1). \end{equation} 
Passing to the
dual  we get the homology isomorphisms as vector spaces
\begin{equation}
 \label{eq:hxa} H_{\ast}(G_{\la}) \cong H_{\ast}(U_0 ) \otimes  H_{\ast}(K_0) \end{equation}
\begin{equation} \label{eq:ha} H_{\ast}(G_{\la}) \cong
H_{\ast}(U_1 ) \otimes H_{\ast}(K_1). \end{equation}

\subsection{The elements $x_i,y_i,t $ and $w_i$}

Denote by $t$ the generator of infinite order in $H_1(
G_{\la};\Z)$, $\la >0$. Recall that $H_{\ast}(SO(3))=\La
(x_1,x_2)$ where $\La$ is the exterior algebra on generators $x_i$
of degree $i$.  Thus $H_{\ast}(K_0)=\La (x_1,x_2,z_1,z_2)$, where
$x_i,z_i$ represent rotation in first and second factors
respectively.  The homology of the $SO(3)$ factor in $K_1 \cong
SO(3) \times S^1$ is generated by $y_i$, and we explain in the next lemma
the relation between  these generators and the generators $x_i$ and $z_i$.
\begin{lemma}
The homology ring of the diagonal in $SO(3) \times SO(3)$, $SO_d(3)$,  is given by
$H_{\ast}(SO_d(3))=\La(y_1,y_2)$ where 
\begin{eqnarray*}
y_1 & = & x_1+z_1 \\
y_2 & = & x_2+z_2+x_1z_1 \\
y_3 & = & x_3+z_3+x_1z_2+x_2z_1,
\end{eqnarray*}
 $x_i$ and $z_i$, with $i=1,2$ are the generators of the homology ring of $SO(3)
\times SO(3) $ and  $x_3=x_1x_2$, $y_3=y_1y_2$ and $z_3=z_1z_2$.
\end{lemma}
\begin{proof}
It is clear that $y_i$ includes terms $x_i+z_i$, just by looking at
the cell structure. Note that  the cup product is defined using the
diagonal map $d\co SO(3) \rightarrow SO(3) \times SO(3) $. If $\al
\in H^{\ast}(SO(3))$ generates $H^1(SO(3))$ then $(\al \cup \al)
(y_2)=d^{\ast}(\al \times \al )(y_2)$. Now we need to define the
duals $\hat{x_1}$ and $\hat{z_1}$ of $x_1$ and $z_1$ respectively.
Let  $\hat{x_1}$ be the element in $H^1(SO(3) \times SO(3))$ such
that $\hat{x_1}(x_1)=1$, $\hat{x_1}(x_i)=0$ if $i=2$ or 3, and
$\hat{x_1}(z_i)=0$ if $i=1,2,3$. $\hat{z_1}$ is defined in a
similar way. We know that the cup product of $\hat{x}_1$ and
$\hat{z}_1$ does
 not vanish, so we have $0 \neq (\hat{x}_1 \cup \hat{z}_1)(y_2)$. Hence
 \begin{eqnarray*}
(\hat{x}_1 \cup \hat{z}_1)(y_2) & = & d^{\ast} (\hat{x}_1 \times
\hat{z}_1)(y_2) \\
& = & (\hat{x}_1 \times \hat{z}_1) (d_{\ast}y_2) \neq 0.
\end{eqnarray*}
Therefore we see that $d_{\ast}y_2$ must have a component in
 $H_1(SO(3))\otimes H_1(SO(3))$. The only element like that is
 $x_1z_1$, so $y_2$ must involve this element.
The result for $y_3$ follows immediately from that for $y_1$ and $y_2$ by multiplication. 
\end{proof}

It follows that the generators $y_i$ commute with generators $x_i$ and
$z_i$. From injections ${i_0}_{\ast}$ and ${i_1}_{\ast}$ we have elements
$t,x_i,z_i$ and $y_i$ in $H_{\ast}(G_{\la} )$. From
isomorphisms \eqref{eq:MV} and \eqref{eq:hxa} we know that the
rank of $H_1(G_{\la} )$ is 3 and as we just showed we have
elements  $t,x_1$ and $y_1$ in $H_1(G_{\la} )$. Clearly these are linearly
independent.

Looking at \eqref{eq:hxa} and \eqref{eq:ha} we see that $t$ must
have a nonzero image in $H_1(U_0)$. On the other hand, since the
homology of the $SO(3)$ factor in $K_1$ is generated by the 
$y_i$, the $x_i$ must have a nonzero image in
$H_i(U_1)$. The class $x_1$ must correspond by \eqref{eq:MV}
to a class in $H_2(U_0)$ and we will see in Lemma \ref{lm:gd}
below that this class is the image  $x_1t$ in $U_0$. $x_1$ is a
spherical representative of the first $SO(3)$ factor in $H_1 ( K_0
)$. Therefore, since $K_0$ acts on $\Jj_{\la}$ by multiplication
on the left there is a well defined 2--cycle $x_1t$ in $U_0$. More
precisely if  $x_1$ is represented by $$ S^1 \rightarrow
 G_{\la}\co u \mapsto x_1(u)$$ and  $t$ by $$S^1 \rightarrow G_{\la}\co v \mapsto
 t(v)$$ we define a 2--cycle in $G_{\la}$ given by the map
 $$S^2=S^1 \times S^1 /S^1 \vee S^1
\rightarrow G_{\la} $$ induced by the commutator $$S^1 \times S^1
\rightarrow G_{\la}\co  (v,u) \mapsto
t(v)x_1(u)t(v)^{-1}x_1(u)^{-1}.$$ It turns out that the projection
of this element in $H_{\ast}(G_{\la})$ to $H_{\ast}(U_0)$, under
the projection map ${p_0}_{\ast}$, is the 2--cycle $x_1t$ in $U_0$.
In order to see that let us recall that for any group
$G$ the Samelson  product   $[x,y] \in \pi_{p+q} (G)$  of elements
$x \in \pi_p(G)$ and $y \in \pi_q(G)$ is represented by the
commutator $$S^{p+q}=S^p \times S^q / S^p \vee S^q \rightarrow G \co 
(u,v) \mapsto x(u)y(v)x(u)^{-1}y(v)^{-1}.$$ The  Samelson product  in
$\pi_{\ast}(G)$  is related to the    Pontryagin product in \linebreak
  $H_{\ast}(G;\Z)$  by the formula  $$[x,y]=xy - (-1)^{|x||y|}yx,$$
where we suppressed the Hurewicz homomorphism $\rho \co  \pi_{\ast}
(G) \rightarrow H_{\ast}(G)$ to simplify the expression. Therefore
we see that this 2--cycle is given by the commutator $[x_1,t]$, so
in homology  is simply given by
$x_1t +tx_1 \in H_2(G_{\la};\Z_2) $. Similarly we define a cycle in $H_4(G_{\la};\Z_2)$
that in homotopy is given
 by the commutator $[t,x_3]$. Although $x_2$ is not a spherical
 class, ie, $ x_2 \notin \pi_2 ( G_{\la})$ we can consider a
 cycle in degree 3 given by $x_2t+tx_2$ in  $H_{\ast}(G_{\la},\Z_2) $.

\begin{definition} \label{def:wi}We define  elements $w_i \in
H_{i+1}(G_{\la},\Z_2) $ to be  the commutators  $x_it+tx_i $ with
$i=1,2,3$. For a word in the $w_i's$ we use the notation
$w_I=w_{i_1}...w_{i_n}$ with $I=(i_1,...,i_n)$. \end{definition}

The reason why we use these classes $x_it+tx_i$ instead of
simply $x_it,tx_i$ is first because they project simultaneously to
additive generators in $H_{\ast}(U_1)$ and $H_{\ast}(U_0)$ so it
is easier to see the correspondence between elements in
isomorphisms \eqref{eq:hxa} and \eqref{eq:ha}. Secondly they are
in the kernel of the subalgebra of  $H^{\ast}(G_{\la}) $ generated
by the duals $\hat{t}$ and $\hat{x_i}$ of $t$ and  $x_i$. We show
this fact in the next lemma, but first we define the duals of these
elements  in $H^1(G_{\la} )$. $\hat{t}$ is the element in
$H^1(G_{\la} )$ such that $\hat{t}(t)=1$ and
$\hat{t}(x_1)=\hat{t}(y_1)=0$. We define $\hat{x_1}$ and
$\hat{y_1}$ in the obvious way. We also have
$\hat{x_i}=(\hat{x_1})^i$ and $\hat{y_i}=(\hat{y_1})^i$.

\begin{lemma} \label{lm:ng} The cup product  $(\hat{t} \cup \hat{x_i})$ evaluated at the commutator $[x_i,t]$ is 0  where $\hat{t}$ and $\hat{x_i}$ represent the dual
of $t$ and $x_i$ in $H^{\ast}(G_{\la}) $ respectively. \end{lemma}
 \begin{proof} Although in this section we are working with $\Z_2$
coefficients we
 will prove a stronger result by showing that the statement is true also
over $\Z$.
  Note that $ (\hat{t} \cup
\hat{x_i}) ( [x_i, t] )= (\hat{t} \cup \hat{x_i}) ( x_it+tx_i)
=(\hat{t} \cup \hat{x_i})( x_it) +(\hat{t} \cup \hat{x_i})(tx_i)$
and we show that $(\hat{t} \cup \hat{x_i})(tx_i)= (\hat{x_i} \cup
\hat{t})(x_it)=- (\hat{t} \cup \hat{x_i})(x_it)=1$. For example,
in the case when $i=1$  consider $f\co  S^1 \times S^1 \rightarrow
G_{\la} \co  (t,s) \mapsto \ph_t \psi_s$, where $S^1\rightarrow
G_{\la}\co t \mapsto \ph_{t}$ and $S^1\rightarrow G_{\la}\co s \mapsto
\psi_{s}$ represent the cycles $t$ and $x_1$ respectively. Then
\begin{eqnarray*} (\hat{t} \cup \hat{x_1})(tx_1) & = &
f^{\ast}(\hat{t} \cup \hat{x_1})[S^1 \times S^1] \\
                               & = & f^{\ast}(\hat{t})\cup
f^{\ast}(\hat{x_1})[S^1 \times S^1] \\
                               & = & f^{\ast}(\hat{t})[S^1]
f^{\ast}(\hat{x_1})[S^1] =1.
\end{eqnarray*} Thus $ (\hat{t} \cup \hat{x_1}) ( [x_1, t] )=
-1+1=0$. \end{proof}

We look at an additive basis for each group $H^k(G_{\la})$ and $H_k(G_{\la})$ in the sense that we want to have a canonical identification between $H^k(G_{\la})$ and $H_k(G_{\la})$, this meaning that if $\{c_{\alpha} \}$ is an additive basis for $H_k(G_{\la})$ then $\{ \hat{c}_{\alpha} \}$ is an additive basis for  $H^k(G_{\la})$ where $\hat{c}_{\alpha}$ is the element satisfying $\hat{c}_{\alpha}(c_{\beta}) = \delta_{\alpha \beta}$.
Using this identification we see from the previous lemma that the dual of the commutators $[x_i, t]$ represent classes in $H^{\ast}(G_{\la})$ which are not in the subalgebra of $H^{\ast}(G_{\la})$ generated by $\hat{t}$, $\hat{x_i}$ and $\hat{y_i}$.

We choose  a normalized set of elements in the subring of
$H_{\ast}(G_{\la})  $ generated by $t,x_i,y_j$ with $i,j=1,2,3$.
\begin{lemma} \label{lm:adb} Any word in the  $t,x_i,y_j$ with
$i,j=1,2,3$ is a sum of elements of the form \begin{equation}
\label{eq:adb} w_It^{\epsilon_t}x_i^{\epsilon_i}y_j^{\eta_j},
\end{equation} where $\epsilon_t$,$\epsilon_i,\eta_j=0$ or $1$,
$I=(i_1,...,i_k)$ and $i,j=1,2 $ or $3$ ($x_3 = x_1x_2, y_3=y_1y_2
)$. \end{lemma} \begin{proof} We know that $y_j$ commutes with all
other elements and we have  equations $$x_iw_j=w_ix_j \; {\rm{if}}
\;(i,j) \neq (1,2) \; {\rm{and}} \;(2,1)$$ $$x_iw_j=w_ix_j +w_3\;
{\rm{if}} \;(i,j) = (1,2) \; {\rm{or}} \;(2,1).$$ We also know
that $x_it = tx_i +w_i$, $t$ commutes with $w_i$ for $i=1,2,3$ and
$t^2=0$. These facts together with the two equations imply that we
can always bring any copy of $x_i$ to the right of $t$ and the ${w_i}'s$,
adding, if necessary, words on the ${w_i}'s$. \end{proof}

\subsection{A generating set for $H_{\ast}(G_{\la})$}

In this subsection the aim is to show that  the elements
$x_i,y_j,t$ generate the ring $H_{\ast}(G_{\la};\Z_2)$. In order
to do that we give a geometric description of the isomorphism
$H_{p+1}(U_0) \cong H_{p}(U_1)$.

 We have projections ${p_i}_{\ast} \co H_{\ast}(G_{\la}) \rightarrow H_{\ast}(U_i)$ with
 $i=0$ or $1$.
 Since $x_i$ has image in $H_i ( K_0 )$ and $t$ has image in $H_1( K_1)$ we can conclude  that ${p_0}_{\ast}([x_i,t])= {p_0}_{\ast}(x_it)$ in $H_{\ast}(U_0)$ and ${p_1}_{\ast}([x_i,t])=
 {p_1}_{\ast}(tx_i)$  in
 $H_{\ast}(U_1)$.
 We write $t$ for ${p_0}_{\ast}(t) \in H_{\ast}(U_0)$ and $x_i$
 for ${p_1}_{\ast}(x_i) \in H_{\ast}(U_1)$. However it will
 be convenient to  distinguish notationally between the
different incarnations of $w_1$,$w_2$,$ w_3$ on the different
spaces. We will denote by $v_i={p_0}_{\ast}(w_i)$ the generators
in $H_{\ast}(U_0)$ and by $u_i={p_1}_{\ast}(w_i)$  the generators
in $H_{\ast}(U_1)$ where $i=1,2$ or 3. Let $v_I={p_0}_{\ast}(w_I)$
and $u_I={p_1}_{\ast}(w_I)$ where $w_I$ is given as in Definition
\ref{def:wi}. This way we give meaning to expressions as
$v_iv_j={p_0}_{\ast}(w_iw_j)$, $v_it={p_0}_{\ast}(w_it)$ and
$u_iu_j={p_1}_{\ast}(w_iw_j)$. We can write $v_it$ or $tv_i$ to
refer to the same element because $t$ commutes with $w_i$ in
$H_{\ast}(G_{\la} )$.
 Note that $H_{\ast}(U_i )$  is a  left $H_{\ast}(G_{\la}
 )$--module, so $H_{\ast}(G_{\la} )$ acts on $H_{\ast}(U_i
 )$ by multiplication on the left. Using this module action we have
  $v_{I'}=w_i.v_I$ and $u_{I'}=w_i.u_I$ for $I'=(i,I)$.

We can choose  right inverses $s_i\co H_{\ast}(U_i ) \rightarrow
H_{\ast}(G_{\la} )$ such that $s_0(t)=t$, $s_0(v_i)=w_i$,
$s_1(x_i)=x_i$, $s_1(u_i)=w_i$  and  ${p_i}_{\ast} \circ s_i=id $.
They exist because of isomorphisms \eqref{eq:hxa} and
\eqref{eq:ha}. Moreover  we can choose $s_i$ such that $s_0$
preserves multiplication by $t,w_i$ and $s_1$ preserves
multiplication by $w_i$.

\begin{lemma} \label{lm:gd} The isomorphism $H_{p+1}(U_0) \cong
H_{p}(U_1)$ is given by the map  $$ \psi \co  H_p (U_1) \rightarrow
H_{p+1}( U_0)\co c \mapsto {p_0}_{\ast}(s_1(c)t) $$ \end{lemma}
 \begin{proof} Note that since $U_1$ is a codimension 2
submanifold of $\Jj_{\la}$, there is a circle bundle $\p
\mathcal{N}_{U_1}$ where $ \p \mathcal{N}_{U_1}$ is a neighborhood
of $U_1$ in $U_0$:
 \begin{equation}
\label{eq:fb}
\xymatrix{ {S^1} \ar[r] & {\partial \Nn_{U_1}}
\ar[d]^{\pi} \\
 & U_1 }
\end{equation}
Therefore for any map representing a cycle $ c \in H_p (U_1)$ we can obtain a cycle in  $ H_{p+1}( U_0)$ by lifting the map to  $\p \mathcal{N}_{U_1}$ using the fibration \eqref{eq:fb}. More precisely, using the section $s_1$ we can lift $c$ to a cycle in $H_*(G_{\la})$. This is represented by a map $\tau\co  C \rightarrow G_{\la}\co  z \mapsto \tau (z) $. Now note that there is a map from the image of $\tau$ to $U_0$  given by $g\mapsto g_{\ast}J$, where $g \in \tau (z)$ and  we can choose $J \in U_0$ close to $U_1$. In fact, we can choose $J \in \mathcal{N}_{U_1}$ so close to $U_1$ such that $g_{\ast}J \in \mathcal{N}_{U_1}$ also. Then using the $S^1$ action, represented by $t$,  we  define a map to $U_0$, representing a cycle in  $\mathcal{N}_{U_1} \subset U_0$: for each $g$ in the image of $\tau$  we   get a loop around $U_1$  defined by
$g_{\ast}(t_{\ast}J)=(gt)_{\ast}J$. Therefore the cycle $c$ lifts to ${p_0}_{\ast}(s_1(c)t)$ which represents a cycle in  $H_*(U_0)$.
\end{proof}

\begin{remark}  \label{rm:gd} {Using the notation introduced above
we can write
\begin{eqnarray*}
\psi(x_i) & = & {p_0}_{\ast}(s_1(x_i)t)=x_it=v_i, \\ \psi (u_i) &
= & {p_0}_{\ast}( s_1(u_i)t)=w_it=v_it, \\ \psi (u_Ix_i) & = &
w_I.v_i=v_{I'} \end{eqnarray*}  with $I'=(I,i)$ and
$\psi(u_I)=v_It$.

The map $\phi$ that gives the corresponding isomorphism  in
cohomology, $$\phi\co H^{p+1} (U_0) \rightarrow H^{p} (U_1),$$ is the
composite of the restriction $ i^{\ast}\co  H^{\ast}(U_0)
\longrightarrow H^{\ast}(\partial \mathcal{N}_{U_1})$ with
integration over the fiber of the projection $\pi \co \partial
\mathcal{N}_{U_1} \rightarrow U_1$ of  the fibration \eqref{eq:fb}
 (see \cite{AM}).}
\end{remark}

Now we use the lemma to prove the following proposition.

\begin{proposition} \label{prop:ring}The generators of the
Pontryagin ring $H_{\ast}(G_{\la}) $ are $t,x_i$, $y_j$ with
$i,j=1,2$.\end{proposition} \begin{proof} The existence of injections ${i_i}_{\ast} \co H_*(K_i) \rightarrow H_{\ast}(G_{\la}) $ with $i=0$ or 1 imply that we have elements $t,x_i,y_j$ in $H_{\ast}(G_{\la}) $. Let $R_{\ast} \subset
H_{\ast}(G_{\la}) $ be the subring generated by $t,x_i,y_j$.
Suppose there is an element of minimal degree $d$ in
$H_{\ast}(G_{\la})  - R_{\ast}$. From isomorphism \eqref{eq:hxa} we can
conclude that such an element would be mapped  to a sum of elements $$\sum_l
c_l \otimes k_l \ \in \
\underset{l}{\oplus}(H_{d-l}(U_0) \otimes H_l(K_0))$$ with
$0\leq l \leq 6$. For some $l$, $c_l$ is not a polynomial in the
$v_I,t$. Take the largest such $l$. By the isomorphism in Lemma
\ref{lm:gd}  and Remark \ref{rm:gd} this would create an
element  in $H_{d-l-1}(U_1)$ that is not a polynomial in $u_I$ and
$x_i$. But this is impossible because this would give rise to a
new generator in $H_{d-l-1}(G_{\la})$ corresponding to this new
element in $H_{d-l-1}(U_1) \otimes H_0(K_1)$ and this contradicts
the minimality of $d$. \end{proof}

\subsection{Main theorem}

We start by showing that we have isomorphisms $H_{\ast}(G_{\la} )
\cong H_{\ast}(U_i) \otimes H_{\ast}(K_i)$ given by the Pontryagin
product.  More precisely, we can define maps \begin{equation}
\label{eq:pdh} \varphi_i \co H_{\ast}(U_i)\otimes H_{\ast}(K_i)
\rightarrow H_{\ast}(G_{\la}) \co c\otimes k \mapsto s_i(c).k
\end{equation} with $i=0$ or 1. Since $K_i \subset G_{\la}$ and $i_i$
is injective in homology we denote ${i_i}_{\ast}(k)$ simply by $k$.
Recall that we have projections $p_i\co G_{\la} \rightarrow U_i$ as
defined in diagram \eqref{eq:diag} and these induce maps
${p_i}_{\ast}\co H_{\ast}(G_{\la}) \rightarrow H_{\ast}(U_i)$ in
homology. It is clear that ${p_i}_{\ast}(s_i(c).k)=0$ if $k \in
H_{\ast}(K_i)$, with $\ast >0$ and the product $s_0(c)k$ is an element
in the normalized set defined in Lemma \ref{lm:adb}, because $s_0(c)$
is a product of $w_i's$ and $t$ and $k$ is a product of $x_i$ and
$y_j$. We now claim that these maps are isomorphisms.

 \begin{proposition}  \label{prop:iso}
 The maps $\varphi_i \co H_{\ast}(U_i)\otimes H_{\ast}(K_i) \rightarrow
H_{\ast}(G_{\la}) \co  c\otimes k \mapsto s_i(c).k$ given by
Pontryagin product are isomorphisms. \end{proposition}

\begin{proof} Consider the elements of the form
$v_It^{\epsilon_t}$, with $\epsilon_t
 =0,1$ in  $H_{\ast}(U_0)$. If they are not linearly independent, choose a
maximal linearly independent subset $B=\{c_\al\}$. It follows from
Proposition \ref{prop:ring} that this is a basis for
$H_{\ast}(U_0)$. Now consider the image in $H_{\ast}(G_{\la}) $ of
$B$. This is given
 by $B'=\{s_0(c_{\al})\}$ with $c_{\al} \in B$. These are  elements of the
form
 $w_It^{\epsilon_t}$, $\epsilon_t =0,1$ and the set $B'$ is
 linearly independent. Therefore it is an additive basis for the space
spanned
 by elements of the form $w_It^{\epsilon_t}$. Note that
$H_{\ast}(G_{\la}) $  has a subalgebra isomorphic to
$H_{\ast}(K_0)$
   and an additive basis for this is
$D=\{k_\ga\}=\{x_i^{\epsilon_i}y_j^{\eta_j}\}$
where $\epsilon_i$
    and $\eta_j$ are equal to 0 or 1, so an additive basis for
$H_{\ast}(G_{\la})$  will contain all elements of this form.
To prove the theorem in the case $i=0$ we need to show that the
set $B''=\{s_0(c_{\al}).k_\ga\}$ where $s_0(c_\al) \in B'$ and $k_\ga
\in D$ is an additive basis of $H_{\ast}(G_{\la})$. We start by proving  that these elements generate additively $H_{\ast}(G_{\la})$.
Suppose we have an element $a \in H_{\ast}(G_{\la})$.
 From Proposition \ref{prop:ring} and Lemma \ref{lm:adb} it is
known  that every element in $H_{\ast}(G_{\la })$ is a sum of
elements of the form \eqref{eq:adb}. Therefore $$a=\sum_\al
w_{J_\al}t^{\epsilon_\al}x_{i_\al}^{\epsilon_i}y_{j_\al}^{\eta_j}.$$
It is also known that $x_{i_\al}^{\epsilon_i}y_{j_\al}^{\eta_j}$ is in $D$
and if $w_{J_\al}t^{\epsilon_\al}$ is not in $B'$ we can write it
as sum of elements in $B'$. Thus $a$ is a sum of elements in
$B''$.

The next step is to show that the elements in $B''$ are linearly
independent. We know that for a fixed degree $d$, the dimension of
$H_d(G_{\la})$ is given by $$\sum_{l=0}^d \dim H_l(U_0)\times \dim
H_{d-l}(K_0),$$ because of the vector space isomorphism
\eqref{eq:hxa}. But this is precisely the number of elements in
$B''$ of degree  $d$. So they must be linearly independent otherwise their
span would not be the space $H_{\ast}(G_{\la })$. This
means that the set $B''=\{s_0(c_{\al}).k_\ga\}$ defined above is
an additive basis for $H_{\ast}(G_{\la})$. Therefore $\varphi_0$
is an isomorphism.

In the case $i=1$, $\varphi_1$ maps $c \otimes k$ to $s_1(c).k$
and this is not an element in the form \eqref{eq:adb}. However we can prove
an analogous result to Lemma \ref{lm:adb} stating that any word in
the $x_i,y_j,t$ is a sum of elements of the form
$w_Ix_{i}^{\epsilon_i}ty_{j}^{\eta_j}$. This is clear because
$w_Itx_{i}y_{j}=w_Ix_ity_j+w_Iw_ity_j$ for all $I,i$ and $j$. 
Repeating the steps for the case $i=0$ and using isomorphism
\eqref{eq:ha} instead of isomorphism  \eqref{eq:hxa} it follows easily that $\varphi_1$ is also an isomorphism.
 \end{proof}

We are now in position to calculate the algebra structure on
$H_{\ast} (G_{\la})$.

 \begin{theorem} \label{tm:main}If $0 <\la
\leq 1$ then $$H_{\ast} (G_{\la} ; \Z_{2} )= \La(y_1,y_2) \otimes
\Z_2\langle t,x_1,x_2\rangle  /R$$ where $\deg y_i= \deg x_i =i$, $\deg t=1$ and $R$ is
the
set of relations $ \{ t^2=x_i^2=0, x_1x_2=x_2x_1 \}$.
\end{theorem} \begin{proof} We already know from Proposition
\ref{prop:ring} that the generators of the Pontryagin ring are
$t,x_i,y_j$. Therefore it is sufficient  to prove that the only relations between
them are the ones in $R$, the commutativity of $y_i$ with  $x_i$ and $t$ plus the ones on $y_i$ coming from the definition of an exterior algebra $\La(y_1,y_2) $. We also know from Lemma \ref{lm:adb}, assuming only these relations, that any word in these generators is  a sum of elements of the form  \eqref{eq:adb}. Thus  if we prove that this set of elements  give an additive basis of  $ H_{\ast}(G_{\la}) $ we prove that there are no more relations between the generators $t,x_i,y_i$, because the existence of another relation  would give rise to one  between the elements of the form \eqref{eq:adb} and  they would not be linearly independent.
We will prove that by induction. The induction hypothesis
is that up to dimension $d-1$ elements of form
  \eqref{eq:adb} are linearly independent, thus there are no relations
between them up to dimension $d-1$.
Suppose there was one of minimal degree $d$ in $
H_d(G_{\la})$. The first step is to show that it would
be between the $w_i's$ only. Assume it was given by a
finite sum of the type $$ \sum_k w_{I_k} A_k = 0 $$ where
$w_{I_k}$ is a word on the $w_i's$ and $A_k=t^{\epsilon_k}b_k$
where $b_k$ is an element in  $H_{\ast}(K_0)$ and $\epsilon_k$
equals $0$ or $1$. Then from Proposition \ref{prop:iso} with $i=0$ it follows  that we must have $$ \sum_k w_{I_k}t^{\epsilon_k}
\otimes b_k= 0 .$$ We can group together the terms in which $b_k$
is the same, thus we can write the relation as $$\sum_k
(\sum_{l_k} w_{l_k}t^{\epsilon_{l_k}})\otimes b_k=0$$ where now
$b_k$ runs over a set of basis elements of $H_{\ast}(K_0)$. This
implies that we have a relation of the type $$ \sum_l
w_{I_l}t^{\epsilon_l}=0.$$ Using Proposition \ref{prop:iso} with $i=1$
we show that it is between the $w_i's$, because $$\sum_l
w_{I_l} \otimes t^{\epsilon_l}=\sum_{l'} w_{I_{l'}} \otimes t+
\sum_{l''} w_{I_{l''}} \otimes 1= 0$$ implies $$\sum_{l'}
w_{I_{l'}}=0 \mbox{ and } \sum_{l''} w_{I_{l''}}=0.$$ A relation
in the $w_i's$ projects, under the map ${p_0}_{\ast}$, to one  on the ${v_i}'s$ in $ H_d(U_0)$.  Using isomorphism
\eqref{eq:MV} this would  give a relation in degree $d-1$ between
the ${u_i}'s$ and $x_i's$ in $H_{d-1}(U_1)$. But this contradicts
the induction hypothesis because such relation implies one in
$H_{\ast}(G_{\la})$ with $\ast$ at most equal to $d-1$.
\end{proof}

\noindent The next corollary is an immediate consequence of the proof of the theorem.

\begin{corollary} The Pontryagin ring $H_{\ast}(G_{\la}) $ contains a free
noncommutative ring on  3 generators, namely $w_1,w_2,w_3$.
\end{corollary}

We proved also the following proposition:

\begin{proposition} An additive basis for $H_{\ast}(G_{\la}) $ is
given by \begin{equation}
w_It^{\epsilon_t}x_i^{\epsilon_i}y_j^{\eta_j},
\end{equation} where $\epsilon_t$,$\epsilon_i,\eta_j=0$ or $1$,
$I=(i_1,...,i_n)$
and $i,j=1,2 $ or $3$ ($x_3 = x_1x_2, y_3=y_1y_2 )$.
\end{proposition}

\subsection{Relation between cohomology and homology}

Establishing the vector space  isomorphisms \eqref{eq:cxa} and
\eqref{eq:ca} does not imply that we have algebra isomorphisms on
cohomology. That is proved in the next lemma.

\begin{lemma} \label{lm:iso}The following isomorphisms hold as
algebra isomorphisms: \begin{equation} \label{eq:iso}
H^{\ast}(G_{\la}) \cong H^{\ast}(U_i)  \otimes H^{\ast}(K_i)
\mbox{ with $i=0,1$ }\end{equation} \end{lemma} \begin{proof} The
proof is based in the argument used by Abreu in \cite{A} with some
necessary changes. $H^{\ast}(G_{\la})$ has subalgebras
${p_i}^{\ast}(H^{\ast}(U_i) ) \cong H^{\ast}(U_i) $. From Theorem
\ref{thm:heq} we know that $\Map_1 (S^2)$ is homotopy equivalent
to $SO(3) \times \Omega$ where $\Omega $ denotes the universal
covering space of ${\Map_1}_{\ast} (S^2)$. Therefore we have a map
$\Map_1 (S^2) \times \Map_1 (S^2) \rightarrow SO(3) \times SO(3)
$. The composite of $G_{\la} \rightarrow \Map_1 (S^2) \times
\Map_1 (S^2)$ with the previous map gives us a map $p\co G_{\ast}
\rightarrow  K_0$. Thus $H^{\ast}(G_{\la})$ has a subalgebra
$p^{\ast}(H^{\ast}(K_0))\cong H^{\ast}(K_0)$. Composing these
inclusions of $H_{\ast}(U_0) $ and $H^{\ast}(K_0)$ as subalgebras
of $H^{\ast}(G_{\la})$ with cup product multiplication in
$H^{\ast}(G_{\la})$ we get a map $$ \nu_0 \co H^{\ast}(U_0) \otimes
H^{\ast}(K_0) \rightarrow H^{\ast}(G_{\la}).$$ $ \nu_0$ is an
algebra homomorphism because $H^{\ast}(G_{\la})$ is commutative
and it is compatible with filtrations  (the obvious one on
$H_{\ast}(U_0) \otimes H^{\ast}(K_0)$ and the filtration $F$ on
$H^{\ast}(G_{\la})$ coming from the fibration on the left in
\eqref{eq:diag}). The degeneration of the spectral sequence at
the $E_2$--term implies that $\nu_0$ is an algebra isomorphism.
This proves  isomorphism \eqref{eq:iso} in the case $i=0$.
 For the case $i=1$ note that the map
$i_1^{\ast}\co H^{\ast}(G_{\la}) \rightarrow H^{\ast}(K_1) $ is
surjective, so there are ${\hat{t}}$ and $\hat{y}$ in
$H^{\ast}(G_{\la})$ such that $i_1^{\ast}(\hat{t})$ and
$i_1^{\ast}(\hat{y})$ generate the ring $H^{\ast}(K_1)$, where
$i_1^{\ast}(\hat{t})$ is the generator of the cohomology of $S^1$
and  $\hat{t}$ is such that $\hat{t}(x_1)=0$.
$i_1^{\ast}(\hat{y})$ is the generator of the cohomology of the
$SO(3)$ factor. Now we need to prove that ${\hat{t}}^2=0$ in
$H^{\ast}(G_{\la})$ in order to claim that the subalgebra of
$H^{\ast}(G_{\la})$ generated by ${\hat{t}}$ and $\hat{y}$ is
isomorphic to $H^{\ast}(K_1)$. \begin{lemma} \label{lm:t2}\qua
${\hat{t}}^2=0$ in $H^{\ast}(G_{\la})$ \end{lemma} \begin{proof}
Using isomorphisms \eqref{eq:MV}, \eqref{eq:hxa} and
\eqref{eq:ha} we can show that the rank of $H_2(G_{\la})$ is 6.
But in $H_2(G_{\la})$ the cycles $x_2,y_2,tx_1,$ $ty_1,x_1y_1,w_1$
are linearly independent. We will show that ${\hat{t}}^2$
evaluated on all these classes is 0. The only one at which is not
obviously 0 is $w_1$. Let the map $\alpha \co S^2=S^1\times S^1/ S^1
\vee S^1 \rightarrow G_{\la} $ represent the 2--cycle $w_1$. Then
${\hat{t}}^2(w_1)=\alpha^{\ast}({\hat{t}}^2)[S^2]=(\alpha^{\ast}({\hat{t}})[
S^2])^2$ and this vanishes because $w_1$ is a spherical class,
ie, $\alpha^{\ast}({\hat{t}}) \in H^1(S^2)=0$. \end{proof} Again
composing these inclusions of $H^{\ast}(K_1)$ and $H_{\ast}(U_1) $
as sub\-algebras of $H^{\ast}(G_{\la})$ with cup product
multiplication we get a map $$ \nu_1 \co H^{\ast}(U_1) \otimes
H^{\ast}(K_1) \rightarrow H^{\ast}(G_{\la})$$ which is an algebra
isomorphism.
 \end{proof}

\begin{remark}
From the isomorphisms in the previous Lemma and in Proposition
\ref{prop:iso} we might be tempted to think that the
diagram \begin{equation} \label{diag:comm} \xymatrix{
H_{\ast}(U_0) \otimes H_{\ast}(K_0) \ar[r]^-{.} \ar[d] &
H_{\ast}(G_{\la}) \ar[d] \\ H^{\ast}(U_0) \otimes H^{\ast}(K_0)
\ar[r]^-{\cup} & H^{\ast}(G_{\la})} \end{equation} commutes, where
the vertical arrows are given by taking a basis $b_i$ of
$H_{\ast}(U_0)$ or $H_{\ast}(G_{\la})$ to its dual basis. Actually
this diagram does not commute. 
To see that we can consider the following example.

We have $$\langle
\hat{w_1} \cup \widehat{x_1y_2},w_1x_1y_2\rangle=1$$ and
 if diagram
\eqref{diag:comm} was commutative then the cup product 
$\hat{w_1} \cup \widehat{x_1y_2}$ evaluated at all other elements
would be  0. But it is not difficult to verify that we  have 
$$\langle \hat{w_1} \cup \widehat{x_1y_2},w_2y_2\rangle  = 1.$$
\end{remark}

\subsubsection{The cohomology ring with $\Z_2$ coefficients}
Since $H^{\ast} (G_{\la} ; \Z_2 )$ is a Hopf algebra we can use
the classification theorem of commutative Hopf algebras over a
field of characteristic 2. It says that $H^{\ast} (G_{\la} ; \Z_2
)$ is a tensor product of the type $$ (\underset{\alpha}{\otimes}
\Z_2[x_{\alpha}]/x_{\alpha}^{h(x_{\alpha})}) \otimes (
\underset{\be}{\otimes} \Z_2[x_{\be}])$$ where $h(x_{\alpha})$ is
a power of 2 (see \cite{MT} for a proof of this structure theorem).

In fact we can prove the following proposition.
\begin{proposition}
If  $0< \la \leq 1$, $H^{\ast}(G_\la;\Z_2)$ is isomorphic, as an algebra, to a tensor product
 of an exterior algebra over $\Z_2$ with an infinite number of generators and a truncated polynomial
 algebra with two generators of multiplicative order 4.
\end{proposition}
\begin{proof}
Note that  Lemma \ref{lm:iso} shows that $H^{\ast} (G_{\la} ; \Z_2 )$ is the tensor product of the algebras $H^{\ast}(U_0)$ and $H^{\ast}(K_0)$. Since $H^{\ast}(SO(3))= \Z_2[\hat{x_1}]/\{{\hat{x_1}}^4=0\}$ we can conclude that $H^{\ast}(K_0)$ is a commutative algebra with two  generators of multiplicative order 4. Next we show that $H^{\ast}(U_0)$ has an infinite number of generators and all have order 2. This  proves the proposition.
As we stated in Remark \ref{rm:gd}, the isomorphism $\phi \co H^{p+1} (U_0) \rightarrow H^{p} (U_1)$ is the
composite of the restriction $$ i^{\ast}\co  H^{\ast}(U_0)
\longrightarrow H^{\ast}(\partial \mathcal{N}_{U_1})$$ with
integration over the fiber of the projection $\pi \co \partial
\mathcal{N}_{U_1} \rightarrow U_1$ of the  fibration
\eqref{eq:fb}.
Since integration over the fiber kills the elements of $\pi^{\ast}(H^{\ast}(U_1))$, for each $0 \neq v \in H^p(U_0)$ there is $k \in \Z_2$ such that  $$i^{\ast}(v)= \hat{t} \cup \pi^{\ast}(u) + k \pi^{\ast}(u')$$  where $u \in H^{p-1}(U_1)$, $u' \in H^p(U_1)$ and $\phi(v) =u$. Therefore we have $ i^{\ast}(v^2)=2\hat{t} \cup \pi^{\ast}(uu') + k^2 \pi^{\ast}({u'}^2)=k^2 \pi^{\ast}({u'}^2) $, because $\hat{t}^2=0$ (proved in Lemma \ref{lm:t2}). Thus $\phi(v^2)=0$. Since $\phi$ is an isomorphism we get $v^2=0$. From knowing that all generators in $H^{\ast}(U_0)$ have multiplicative order 2 it follows that
$H^{\ast}(U_0)$ must have an infinite number of generators, just by comparing the dimensions $\dim H^p(U_0)=\dim H_p(U_0)$ for each $p$. Recall that  Theorem \ref{tm:main}  implies that $H_{\ast} (G_{\la} ; \Z_2 )$ contains a free noncommutative ring on 3 generators that projects to non-zero elements in $H_{\ast}(U_0, \Z_2)$ and the dimensions $\dim H_p(U_0)$ increase as $p$ increases.
\end{proof}

\begin{remark} {In the rational case Abreu proved in \cite{A} that the cohomology ring $H^{\ast}(U_0, \Q)$ contains a generator in dimension 4 of infinite order.  The previous result does not contradict this fact, but shows that  $H_{\ast} (G_{\la} ; \Z )$ contains a divided polynomial algebra. We will confirm this fact in Section 3 when we compute the homotopy type of $G_{\la}$.}
\end{remark}

\section{Homotopy type of $G_{\la}$}

In this section we will show that $G_{\la}$ is homotopy equivalent to the product $X= L \times S^1 \times SO(3) \times SO(3)$ where $L=\Omega \Sigma (S^1 \times SO(3))$. We start by giving some considerations about torsion in $H_{\ast} (G_{\la} ; \Z )$. In the second subsection we explain why we consider the loop space $L$ and the last subsection is devoted to the construction of the homotopy equivalence between $G_{\la}$ and $X$.

\subsection{Torsion in $H_{\ast} (G_{\la} ; \Z )$}\label{sec:torsion}
We can repeat the argument used in Section 2 to compute the Pontryagin rings  $H_{\ast} (G_{\la} ; \Q )$ and $H_{\ast} (G_{\la} ; \Z_p )$, with $\Q$ and $\Z_p$ coefficients, with $p$ prime and $\neq 2$. In this case the homology of $SO(3)$ is given by a single generator in dimension 3. Therefore it is easy to see  that the  generators of the homology ring of $G_{\la}$, in this case,  are simply $t,x_3$ and $y_3$, where $t^2=x_3^2=y_3^2=0$, $y_3$ commutes with $x_3$ and $t$, but $t$ does not commute with $x_3$. Thus we obtain a theorem analogous to Theorem \ref{tm:main}:
\begin{theorem}
If $0< \la \leq 1$ then $$H_{\ast} (G_{\la} ; F ) \cong \La(y_3) \otimes F \langle t,x_3 \rangle /R $$
where $\deg y_3 =3$, $R$ is the set of relations $\{ t^2=x_3^2=0 \}$ and $F$ is the field $\Q$ or $\Z_p$ with $p \neq 2$.
\end{theorem}

\noindent Moreover an additive basis for the homology is given by elements of the form $$w_3^kt^{\epsilon_t}x_3^{\epsilon}y_3^{\eta},$$ where $\epsilon_t,\epsilon, \eta =0$ or 1, $k \in \N$ and $w_3$ is obtained as the commutator of $x_3$ and $t$. Since the results  are  the same if we consider $\Q$ coefficients  or $Z_p$ coefficients with $p \neq 2$  we can conclude that $H_{\ast} (G_{\la} ; \Z )$ has no $p$--torsion if $p \neq 2 $.

 \subsection{ The James construction}
\label{sec-james}
We proved in Section 2 that $H_{\ast} (G_{\la} ; \Z_2 )$ contains a free non-commutative algebra on 3 generators in dimensions 2, 3 and 4. Now the aim is to find an $H$--space $L$ such that the homology ring $H_{\ast}(L;\Z_2)$ is isomorphic to this algebra $\Z_2 \langle w_1,w_2,w_3 \rangle  $ where $w_i$ is in dimension $i+1$. In order to find such space we will use the James construction that we describe next.

For a pointed topological space $(X,*)$, let
$J_k (X) = X^k/{\sim}$ where $$ (x_1,...,x_{j-1}, \ast, x_{j+1},
..., x_k) \sim (x_1,...,x_{j-1},  x_{j+1} ,\ast,..., x_k).$$ The
James construction on $X$, denoted $J(X)$ is defined by $$J(X) =
\lim_{ \stackrel{\rightarrow}{k} }J_k(X),$$ where $J_k(X) \subset
J_{k+1}(X)$ by adding $*$ in the last component. There is a canonical
inclusion $X= J_1(X) \hookrightarrow J(X)$. $J(X)$ is a
topological monoid and any map from $X$ to a topological monoid
$M$ extends uniquely to a morphism $J(X) \rightarrow M$ of
topological monoids. That is, $X \hookrightarrow J(X)$ is
universal with respect to maps from $X$ to topological monoids,
ie, if $f\co X \rightarrow M$ is given there is a unique $
\tilde{f}$  such that the following diagram commutes:$$\xymatrix { X \ar [d]  \ar [dr]^f
\\ JX \ar [r]^{\tilde{f}} & M }$$
$\tilde{f}$ is defined by
$\tilde{f}(x_1,...,x_k)=f(x_1) \ ...\ f(x_k)$.  From the definition, \linebreak $
J^kX/J^{k-1}X = X\wedge X \wedge ... \wedge X $ and  since we have
a filtration $$ JX \supset ... \supset J^kX \supset J^{k-1}X
\supset ...$$ applying the K\"{u}nneth Theorem we conclude the
following (see \cite{W}): 

\begin{theorem} \label{tm:J1}$$ H_{\ast} ( J^kX; {\Z}_2 )= H_{\ast}
( J^{k-1}X; {\Z}_2 ) \oplus H_{\ast}  ( X\wedge X \wedge ...
\wedge X ; {\Z}_2 ) $$ and $$ \tilde{H}_{\ast} (JX; {\Z}_2 )=
{\oplus}_k \left( \tilde{H}_{\ast} (X; {\Z}_2 ) \right)^{\otimes_k} = T (\tilde{H}_{\ast} (X; {\Z}_2
))$$ where given a vector space $H$, $T(H)$ is the tensor
 algebra on $H$ and the last isomorphism is an
isomorphism of Pontryagin rings. \end{theorem} Now note  the
following theorem (see proof in \cite{S}).
\begin{theorem}[James] \label{tm:J2}If $X$ has the homotopy type of a connected
CW--complex then $JX$ and $\Omega  \Sigma X$ are homotopy
equivalent. \end{theorem} The Theorems  \ref{tm:J1} and \ref{tm:J2} imply   that
$$\tilde{H}_{\ast} (\Omega \Sigma X; {\Z}_2 ) \cong T
(\tilde{H}_{\ast} (X; {\Z}_2 ))$$ so, in particular, if $X =S^1
\wedge SO(3)$  we get $$\tilde{H}_{\ast} (\Omega \Sigma (S^1
\wedge SO(3)) ; {\Z}_2 ) \cong T (\tilde{H}_{\ast} ( S^1 \wedge
SO(3); {\Z}_2 )) \cong {\Z}_2 \langle w_1, w_2,w_3 \rangle , $$
where $ w_1, w_2,w_3$ are generators in dimension $ 2, 3,4 $
respectively. So we see that the homology with ${\Z}_2$
coefficients of this space is isomorphic to a subalgebra of
$H_{\ast} (G_{\la} ; {\Z}_2 )$.

 \subsection{The homotopy equivalence}
\label{sec:homtype} The following result is well known: see \cite{Ha}, 
Corollary 3.37.
\begin{proposition} \label{prop:hom} A map $f\co X \rightarrow Y $
induces isomorphisms on homology with $\Z$  coefficients iff it
induces isomorphisms on homology with $\Q$ and $\Z_p$ coefficients
for all primes $p$.
\end{proposition} We now  define the map  $f$ from $X=\Omega \Sigma (S^1 \times SO(3)) \times S^1 \times SO(3) \times SO(3) $ to
$G_{\lambda}$  that  induces isomorphisms on homology with $\Q$
and $\Z_p$ coefficients, for all primes $p$. We have an inclusion
map $$i\co S^1 \times SO(3) \rightarrow G_{\lambda}$$ given by
\begin{equation} \label{eq:samelson}
 (x,y) \mapsto i_1(x)i_0(y)i_1(x)^{-1}i_0(y)^{-1}
\end{equation} where $i_0$ and $i_1$ are the inclusions defined  in
Section \ref{sec:ring}. More precisely, in this formula $i_1$ is
the restriction of the inclusion $K_1 \hookrightarrow G_{\la}$ to
the $S^1$ factor and $i_0$ is the restriction of the inclusion
$K_0 \hookrightarrow G_{\la}$ to the first $SO(3)$ factor. The
restriction to $S^1 \vee SO(3)$ of $i$ is the identity so there is
an  induced map $$ h\co S^1 \wedge SO(3) \rightarrow G_{\lambda}.$$
This map induces the right correspondence between generators in
homology $$h_{\ast}\co  H_{\ast}(S^1 \wedge SO(3);\Z_2) \rightarrow
H_{\ast}(G_{\lambda};\Z_2),$$ this meaning that the three
generators of  $H_{\ast}(S^1 \wedge SO(3);\Z_2)$ are mapped to  $
w_1,w_2,w_3 \in H_{\ast}(G_{\lambda};\Z_2)$, because as we saw
before these generators in $ H_{\ast}(G_{\lambda};\Z_2)$ are
obtained as commutators  of the form  \eqref{eq:samelson}. Moreover
there is a unique map $\tilde{h}$ that extends $h$ to $\Omega
\Sigma (S^1 \wedge SO(3))$ as we explained in Section \ref{sec-james}.
Therefore the map
 \begin{equation} \label{eq:exth}
\tilde{h}_{\ast}\co 
H_{\ast}(\Omega \Sigma (S^1 \wedge SO(3));\Z_2) \rightarrow
H_{\ast}(G_{\lambda};\Z_2)
\end{equation}
takes the generators of $H_{\ast}(\Omega \Sigma (S^1 \wedge SO(3));\Z_2)$ to
the elements $w_1,w_2,w_3$ in  $H_{\ast}(G_{\lambda};\Z_2)$.
Now consider the map $f\co L \times S^1\times
SO(3)\times SO(3) \rightarrow G_{\lambda}$ given by 
\begin{equation} \label{eq:deff}
(w,t,y,x) \mapsto \tilde{h}(w)i_1(t,y)i_0(x),
\end{equation}
 where $w \in L=\Omega \Sigma (S^1 \wedge SO(3))$.

\begin{lemma}
The map $f$ defined above induces isomorphisms on homology with $\Q$ and
$\Z_p$ coefficients for all primes $p$.
\end{lemma}
\begin{proof}
The map $f$ restricted to $S^1\times SO(3)$ or the second $ SO(3)$ factor is just the
inclusion in $G_{\lambda}$. Moreover, using the K\"{u}nneth formula for
homology with coefficients in a field F, we get 
\begin{equation} \label{eq:Kf}
H_n(X;F)\cong  \underset{p+q+l=n}{\oplus} H_p(L;F) \otimes H_q(S^1\times
SO(3);F) \otimes H_l( SO(3);F).
\end{equation}
Let $F$ be $\Q$ or $\Z_p$ with $p\neq 2$.  Note that in this case
$H_{\ast}(SO(3); F)$ has only a generator in dimension 3. Therefore
 an additive  basis for $H_{\ast}(G_{\lambda};F)$ is given by
$$w_3^kt^{\epsilon_t} x_3^{\epsilon} y_3^{\eta}$$ where $\epsilon_t,
\epsilon, \eta = 0$ or $1$.  Thus  comparing equation \eqref{eq:Kf} and an
additive basis for $H_{\ast}(G_{\lambda};F)$ we conclude that the homology groups
of $X$ and $G_{\la}$ are the same. We just need to show that $f$ induces
those isomorphisms. The elements  $t$, $x_3$ and $y_3$ are the images
of the generators of $H_{\ast}(S^1\times SO(3);F)$ and $H_{\ast}(
SO(3);F)$ under the injective maps
$${i_1}_{\ast}\co H_{\ast}(S^1\times SO(3);F)\rightarrow H_{\ast}(
G_{\lambda};F)$$ and $${i_0}_{\ast}\co H_{\ast}(
SO(3);F)\rightarrow H_{\ast}( G_{\lambda};F)$$ induced by
inclusions $i_0$ and $i_1$. On the other hand the restriction of $f$ to $L$
is given by the map $\tilde{h}$ and we know that
$\tilde{h}_{\ast}$      maps the 4--dimensional generator of $$
H_{\ast}(\Omega \Sigma (S^1 \wedge SO(3));F) \cong
T(\tilde{H}_{\ast}(S^1 \wedge SO(3);F)) \cong F[w_3]$$ to the
element in  $H_4(G_{\lambda};F)$ obtained as the Samelson product of $t$ and
$x_3$. This proves
that $f$ induces an isomorphism in homology with $\Q$ and $\Z_p$
coefficients for all primes $p$ with $p \neq 2$.

 If $F= \Z_2$ then
an additive basis for $H_{\ast}(G_{\lambda};\Z_2)$ is given by the set of elements of the form
\eqref{eq:adb}. It follows from equation \eqref{eq:Kf}  that the
homology groups $H_{\ast}(G_{\lambda};\Z_2)$ and $H_{\ast}(X;\Z_2)$ are
isomorphic. The elements  $t,y_1,y_2$ are  the images of the generators of
$H_{\ast}(S^1 \times SO(3);\Z_2)$ and $x_1,x_2$ are the images of the
generators of $H_{\ast}(SO(3);\Z_2)$. In this case the map
$\tilde{h}_{\ast}$ stated in \eqref{eq:exth} takes the generators of $
H_{\ast}(\Omega \Sigma (S^1 \wedge SO(3));\Z_2) $ to the elements
$w_1,w_2,w_3$ which are  the three generators of the free
noncommutative subalgebra of  $H_{\ast}(G_{\lambda};\Z_2)$.
Therefore we get  another isomorphism in homology.
\end{proof}
Since the conditions of Proposition \ref{prop:hom} are satisfied  we can conclude that the map $f$ defined in \eqref{eq:deff} induces
isomorphisms on homology with $\Z$ coefficients. Using the fact the  $G_{\lambda}$ and $X$ are  both $H$--spaces  it follows that $\pi_1$ acts trivially on all ${\pi_n}'s$ in each one of the spaces (see \cite{W} , pp 119). This allow us to apply the following corollary of Whitehead's theorem ({\cite{Ha} proposition 4.48}): 
\begin{corollary} If $X$ and $Y$ are abelian
CW--complexes   (i. e. $\pi_1$ acts trivially on all
${\pi_n}'s$) then a map $f\co X \rightarrow Y$ that induces
isomorphisms in homology is a homotopy equivalence.
\end{corollary} Therefore we have proved our main  theorem:
\begin{theorem} If $0<\la \leq 1$, $G_{\lambda}$ is homotopy
equivalent to the  product $\Omega \Sigma (S^1 \wedge SO(3))
\times S^1 \times SO(3) \times SO(3) $. \end{theorem}
\begin{remark} {{Although the spaces $G_{\lambda}$ and  $
\Omega \Sigma (S^1 \wedge SO(3)) \times S^1 \times SO(3) \times
SO(3) $ are homotopy equivalent the above homotopy equivalence is
not an $H$--map, because it does not preserve the product
structure. This can be seen by comparing the Pontryagin products
on integral homology. It would be an interesting question to find
an easily understandable $H$--space with a  Pontryagin ring
isomorphic to the one of $G_{\la}$.}} \end{remark}

\end{document}